\documentclass[10pt]{article}
\usepackage{geometry}                % See geometry.pdf to learn the layout options. There are lots.
\usepackage{graphicx}
\usepackage{amssymb,amsthm}
\usepackage{epstopdf}
\usepackage{pdfsync}
\usepackage{graphicx, color}

\newtheorem{definition}{Definition}[section]
\newtheorem{lemma}[definition]{Lemma}
\newtheorem{theorem}[definition]{Theorem}
\newtheorem{proposition}[definition]{Proposition}

\newtheorem{remark}[definition]{Remark}
\newtheorem{example}[definition]{Example}

\def\Proof{\noindent{\it Proof.}\ }

\def\e{\varepsilon}
\def\dx{\,dx}

\def\ZZ{\mathbb{Z}}
\def\NN{\mathbb{N}}
\def\rr{\mathbb{R}}
\def\HH{\mathcal{H}}
\def\dHH{{\rm d}\mathcal{H}^{d-1}}

\def\loc{{\rm loc}}
\def\hom{{\rm hom}}
\def \trait (#1) (#2) (#3){\vrule width #1pt height #2pt depth #3pt}
\def \qed{\hfill
        \trait (0.1) (6) (0)
        \trait (6) (0.1) (0)
        \kern-6pt
        \trait (6) (6) (-5.9)
        \trait (0.1) (6) (0)
\medskip}

\begin{document}
\title{Discrete double-porosity models for spin systems}
\author{Andrea Braides \\Dipartimento di Matematica, Universit\`a di Roma Tor Vergata
\\ via della ricerca scientifica 1, 00133 Roma, Italy\\ \\ Valeria Chiad\`o Piat\\ Dipartimento di Matematica, Politecnico di Torino \\  corso Duca degli Abruzzi 24, 10129 Torino, Italy
\\ \\ Margherita Solci\\
DADU, Universit\`a di Sassari\\
 piazza Duomo 6, 07041 Alghero (SS), Italy
}
\date{}                                           % Activate to display a given date or no date

\maketitle

\abstract{\noindent We consider spin systems between a finite number $N$ of ``species''
or ``phases'' partitioning a cubic lattice $\mathbb{Z}^d$. We suppose that 
interactions between points of the same phase are coercive, while between point of  different phases (or, possibly, between points of an additional ``weak phase'') are
of lower order. Following a discrete-to-continuum approach we characterize the limit as a continuum energy defined on $N$-tuples of sets (corresponding to the $N$ strong phases) composed of a surface part, taking into account homogenization at the interface of each strong phase, and a bulk part which describes the combined effect of lower-order terms, weak interactions
between phases, and possible oscillations in the weak phase.


}

\section{Introduction} 

\begin{figure}[h!]
\centerline{\includegraphics [width=2.5in]{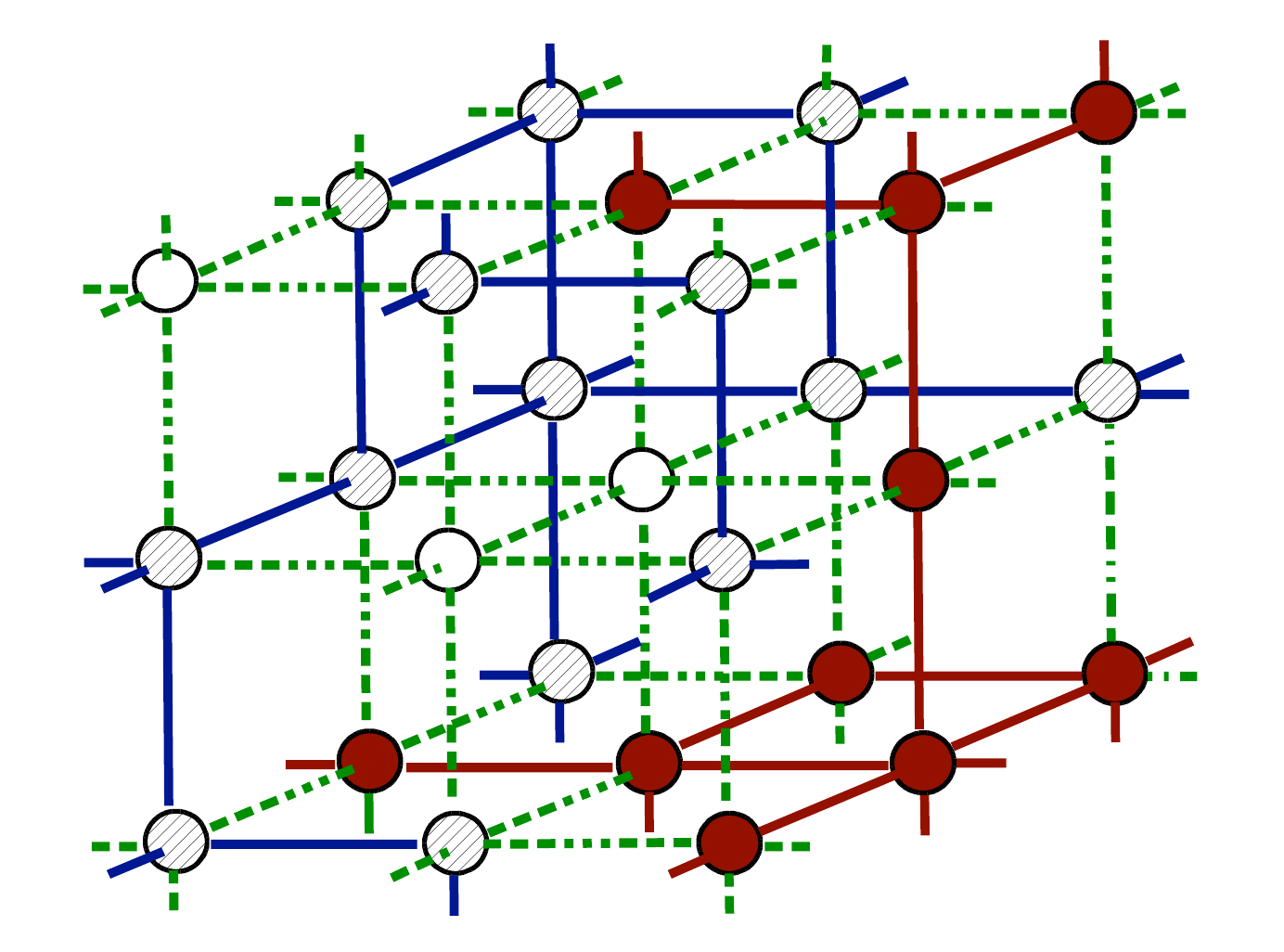}}
\caption{picture of a double-porosity system}\label{doublep}
   \end{figure}

\begin{figure}[h!]
\centerline{\includegraphics [width=2.5in]{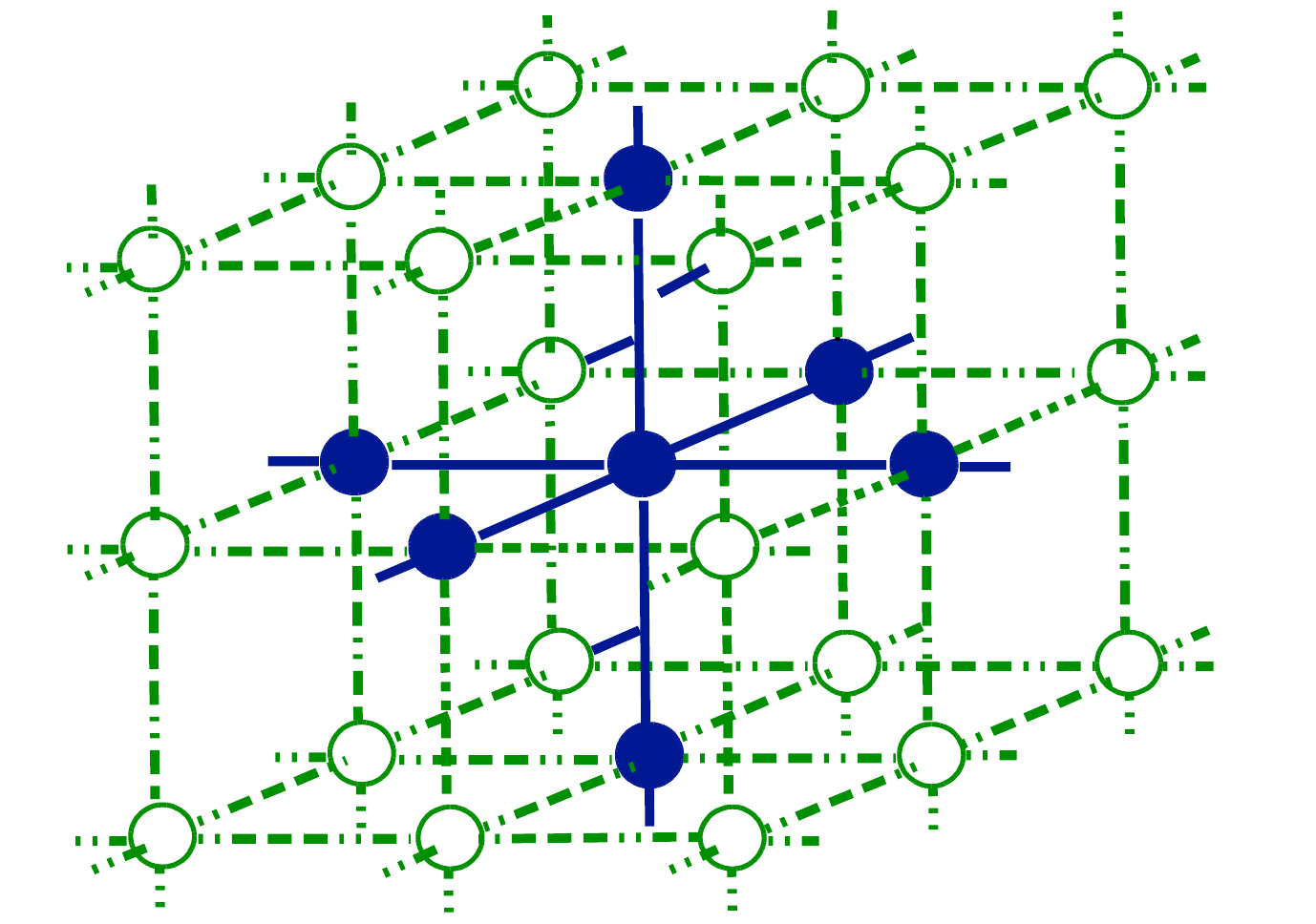}}
\caption{a double-porosity system with an infinite connected weak component}\label{doublep1}
   \end{figure}

In this paper we consider  lattice spin energies mixing  strong ferromagnetic interactions and weak (possibly, antiferromagnetic) pair interactions.  
The geometry that we have in mind is a periodic system of interactions, as that whose periodicity cell is represented in Fig.~\ref{doublep}.
In that picture, the strong interactions between nodes of the lattice (circles)
are represented by solid lines and weak ones by dashed lines.
In this particular case, we have two three-periodic systems of `strong sites'; 
i.e., sites connected by strong interactions, and isolated `weak sites' 
(pictured as white circles). 
   Note that we may also have 
one or more infinite systems of connected weak interactions
as in Fig.~\ref{doublep1}.
In a discrete environment the topological
requirements governing the interactions between the strong and
weak phases characteristic of continuum high-contrast models
are substituted by assumptions on long-range interactions. 
In particular, contrary to the continuum case, for discrete systems with
second-neighbour (or longer-range) interactions we may have a limit
multi-phase system even in dimension one (see the examples in the final section).

This paper is part of a general study of spin systems by means of variational techniques
through the computation of continuum approximate energies,
for which homogenization results have been proved in the ferromagnetic case (i.e, 
when all interactions are strong) by Caffarelli and de la Llave \cite{CDL} and Braides 
and Piatnitski \cite{BP1}, and a general discrete-to-continuum theory of representation and optimization
has been elaborated (see 
%the forthcoming book \cite{ABC} and 
the survey article
\cite{ICM}). In particular, a discrete-to-continuum compactness result and an integral representation of the limit
by means of surface energies defined on sets of finite perimeter has been proved 
by Alicandro and Gelli \cite{AGe}. In that result, the coercivity of energies is obtained by
assuming that nearest neighbours are always connected through a chain of strong interactions.
Double-porosity systems can be interpreted as energies for which this condition does not hold,
but is satisfied separately on (finitely many) infinite connected components. 

Double-porosity problems had been previously considered on the continuum for integral energies (see e.g.~\cite{ADH,BLM,BCP,Pa91,PP,Sa97}) and for interfacial energies (see the works by Solci \cite{Solci,Solci2}). 
A study of discrete double-porosity models 
in the case of elastic energies has been recently carried on in \cite{BCP2}. With respect to
that paper we remark that the case of spin systems allows a very easy proof of an extension lemma
from connected discrete sets, and at the same time permits to highlight the possibility to include a weak
phase with antiferromagnetic interactions, optimized by microscopic oscillations.

\medskip

We are going to consider energies defined on functions parameterized on the cubic lattice $\ZZ^d$ 
of the following  form
\begin{equation}\label{trre}
F_\e(u)=\hskip-.5cm
\sum_{(\alpha,\beta)\in\e {\cal
N}_1\cap(\Omega\times\Omega)}\e^{d-1}a^\e_{\alpha\beta}{(u_\alpha-u_\beta)^2}
+\hskip-.5cm\sum_{(\alpha,\beta)\in\e {\cal
N}_0\cap(\Omega\times\Omega)}\e^d
a^\e_{\alpha\beta}{(u_\alpha-u_\beta)^2} +\sum_{\alpha\in
\Omega\cap\e\ZZ^d}\e^d g(u_\alpha),
\end{equation}
where $\Omega$ is a regular open subset of $\rr^d$,
and $u_\alpha\in\{-1,+1\}$ denote the values of a spin function.
For explicatory purposes, in this formula and the rest of the Introduction, we use a simplified notation
with respect to the rest of the paper, defining
$u=\{u_\alpha\}$  on the nodes of $\Omega\cap\e\ZZ^d$ (instead of, equivalently,
on the nodes of ${1\over\e}\Omega\cap\ZZ^d$).
We denote by ${\cal N}_1$ the set 
of pairs of nodes in $\ZZ^d\times \ZZ^d$
between which we have strong interactions, and by ${\cal N}_0$ the
set of pairs in $\ZZ^d\times \ZZ^d$ between which we have weak
interactions; the difference between these two types of interactions 
in the energy is the scaling factor: $\e^{d-1}$ for strong interactions
and $\e^{d}$ for weak interaction. We suppose that all coefficients
are obtained by scaling fixed coefficients on $\ZZ^d$; i.e.,
\begin{equation}\label{ch-c}
a^\e_{\alpha \beta}= a_{\alpha/\e\,\beta/\e}
\qquad\hbox{ if }\alpha,\beta\in \e\ZZ^d,
\end{equation}
and $a_{jk}$ are periodic of some integer period $T$.
Moreover, we assume that the coefficients of the strong interactions
are strictly positive; i.e., $a_{jk}>0$ if $(j,k)\in {\cal N}_1$.
The `forcing' term containing $g$
and depending only on the point values $u_\alpha$
is of lower-order with respect of strong interactions, but of the
same order of the weak interactions.

We suppose that there there are $N$ infinite connected components
of the graph of points linked by strong interactions, which we
denote by $C_1,\ldots, C_N$. Note that
weak interactions in ${\cal N}_0$ are due either to the existence
of ``weak sites'' or to weak bonds between different ``strong
components'', and, if we have more than one strong graph, the
interactions in ${\cal N}_0$ are present also in the absence of a
weak component.

If we consider only the strong interactions restricted to each strong connected component $C_j$,
we obtain energies 
\begin{equation}\label{trrej}
F^j_\e(u)=
\sum_{(\alpha,\beta)\in\e {\cal
N}^j_1\cap(\Omega\times\Omega)}\e^{d-1}a^\e_{\alpha\beta}{(u_\alpha-u_\beta)^2}
\end{equation}
where ${\cal N}^j_1$ is the restriction to $C_j\times C_j$ of the set ${\cal N}_1$.
This is a discrete analog of an energy on a {\em perforated domain},
the perforation being $\ZZ^d\setminus C_j$. 

We prove an extension lemma that allows to define 
for each $j\in\{1,\ldots,N\}$ a {\em discrete-to-continuum convergence} of 
(the restriction to $C_j$ of) a sequence of function $u^\e$
to a function $u^j\in BV(\Omega;\{\pm1\})$, which is 
compact under a equi-boundedness assumptions for
the energies $F^j_\e(u^\e)$. Thanks to this lemma, such energies
behave as ferromagnetic energies with positive coefficients on 
the whole $\ZZ^d$, which can be {\em homogenized} thanks to
\cite{BP1}; i.e., their $\Gamma$-limit with respect to the convergence
$u^\e\to u^j$ exists, and is of the form
\begin{equation}\label{trre0j}
F^j(u^j)=\int_{S(u^j)\cap\Omega} f^j_{\rm hom}(\nu_{u^j})d\HH^{d-1},
\end{equation}
where $S(u^j)$ is the set of jump points of $u^j$, that can also be
interpreted as the interface between $\{u^j=1\}$ and $\{u^j=-1\}$.

Taking separately into account the restrictions of $u^\e$ to
all of the components $C_j$, we define a {\em vector-valued} 
limit function $u=(u^1,\ldots, u^N)$ and a convergence
$u^\e\to u$, and consider the $\Gamma$-limit with respect
of the whole energy with respect to that convergence.
The combination of the weak interactions and the forcing term give rise to
a term of the form
$$
\int_\Omega \varphi(u)\dx
$$
depending on the values of all components of $u$.
In the case that $\bigcup_{j=1}^N C_j$ is the whole $\ZZ^d$ ,the 
function $\varphi(z^1,\ldots,z^N)$
is simply computed as the average of the $T$-periodic function
$$
i\mapsto\sum_{k\in\ZZ^d} a_{ik}{(u_i-u_k)^2} +g(u_i)
$$
where $u$ takes the value $z^j$ on $C_j$. Note that with this conditions only
(weak) interactions between different $C_j$ are taken into account.
Note moreover that the restriction of the last term $g$ to $\e C_j$ 
is continuously converging to
$$
K_j\int_\Omega g(u^j)\dx,
$$
where $K_j=T^{-d}\#\{i\in C^j; i\in\{0,\ldots, T\}^d\}$ is the 
percentage of sites in $C_j$. 
In general, $\varphi$ is obtained by optimizing the combined effect of
weak pair-interactions and $g$ on the free sites in the
complement of all $C_j$.

Such different interactions can be summed up to describe the 
$\Gamma$-limit of $F_\e$ that finally takes the form
\begin{equation}\label{trre0hom}
F_{\rm hom}(u)=\int_{S(u)\cap\Omega} f_{\rm hom}(\nu_{u})d\HH^{d-1}+
\int_\Omega\varphi(u)\dx,
\end{equation}
where  $f_{\rm hom}(\nu)=\sum_{j=1}^N f^j_{\rm hom}(\nu)$.

We note that the presence of two terms of different dimensions
in the limit highlights the combination of bulk homogenization effects 
due to periodic oscillations besides the optimization of the interfacial structure.
The effect of those oscillations on the variational motions of such systems 
(in the sense of \cite{AGS,2013LN}) is addressed in
\cite{Bso4}.

\section{Notation}\label{Not}
The numbers $d$, $m$, $T$ and $N$ are positive integers.
We introduce a $T$ periodic {\em label function} $J:\ZZ^d\to\{0,1\ldots, N\}$, and the corresponding sets of sites
$$
A_j=\{ k\in\ZZ^d: J(k)=j\},\qquad j=0,\ldots, N.
$$

%
%{\color{red} Eventualmente modificare o notare che $P^0$ non \`e definita solo per $A_0$: tenere conto che i punti forti possono essere connessi deboli alle altre componenti - forse un esempio di perch\'e P dipende da k}
Sites interact through possibly long (but finite)-range interactions, whose range is defined through 
a system $P^j=\{P_k^j\}$ of finite subsets
$P_k^j\subset \ZZ^d$, for $j=0,\ldots, N$ and $k\in A_j$. 
We suppose that

$\bullet$ 
({\em $T$-periodicity}) $P_{k+m}^j=P_k^j$ for  all $m\in T\ZZ^d$;

$\bullet$  ({\em symmetry}) if $k\in A_j$ for $j=1,\ldots N$ ({\em hard components})
and $i\in P_k^j$ then $k+i\in A_j$ and $-i\in P_{k+i}^j$, 
and that $0\in P_k^j$. 

%If $j'\in\{1,\ldots N\}$ and $j'\neq j$ then $P_k^{j'}=\emptyset$.

We say that two points $k,k'\in A_j$ are $P^j$-{\em connected} in $A_j$ if there exists a path $\{k_n\}_{n=0,\ldots,K}$
such that $k_n\in A_j$, $k_0=k$, $k_K=k'$ and $k_n-k_{n-1}\in P_{k_{n-1}}^j$.

We suppose that 

$\bullet$ ({\em connectedness}) 
there exists a unique infinite $P^j$-connected component of each $A_j$ for $j=1,\ldots, N$, 
which we denote by $C_j$. 

Clearly, the connectedness assumption is not a modelling restriction upon introducing more labelling parameters, if the number of infinite connected components is finite. 
Note that we do not make any assumption on $A_0$ and $P^0$. In particular, if 
$k\in A_j$ for $j=0,\ldots N$ and $i\in P_k^0$ then $k+i$ may belong to any $A_{j'}$ 
with $j'\neq j$.

We consider the following sets of bonds between sites in $\ZZ^d$:
for $j=1,\ldots, N$
$$
N_j=\{(k,k'): k,k'\in A_j,  %C_j,
k'-k\in P_k^j\setminus \{0\}\};
$$
for $j=0$
$$
N_0=\{(k,k'): k'-k\in P^0_k\setminus\{0\}, J(k)J(k')=0\hbox{ or } J(k)\neq J(k')\}.
$$
Note that the set $N_0$ takes into account interactions not only among points of the set $A_0$, but also among pair of points in different $A_j$. A more refined notation could be introduced by defining range of interactions $P^{ij}$ and the corresponding sets $N_{ij}$, in which case the sets $N_j$ would correspond to $N_{jj}$ for $j=1,\ldots, N$ and $N_0$ the union of the remaining sets. However, for simplicity of presentation we limit our notation to a single index.

We consider interaction energy densities associated to positive
numbers $a_{kk'}$ for $k,k'\in\ZZ^d$, and forcing term $g$. 
We suppose that  for all $k,k'\in\ZZ^d$

\smallskip
(i) ({\em coerciveness on the hard phase}) there exists $c>0$ such that $a_{kk'}\ge c>0$
if $k\in C_j$ and $k'-k\in P_k^j$ for $j\ge 1$;

(ii) ({\em $T$-periodicity}) $a_{k+m\,k'+m}=a_{kk'}$  for all $m\in T\ZZ^d$;

(iii)  ({\em symmetry})  $a_{k'k}=a_{kk'}$;

(iv) ({\em $T$-periodicity of the forcing term}) $g(k+m,1)= g(k,1)$ 
and $g(k+m,-1)=g(k,-1)$ for all $m\in T\ZZ^d$.

Note that we do not suppose that $a_{kk'}$ be positive for weak interactions.
They can as well be negative, thus favouring oscillations in the weak phase.

\bigskip
Given $\Omega$ a bounded regular open subset of $\rr^d$, 
for $u:{1\over\e}\Omega\cap\ZZ^d\to\{+1,-1\}$ we define
the energies 
\begin{eqnarray}\nonumber\label{fep}
F_\e(u)=F_\e\bigl(u,{1\over\e}\Omega\bigr)&=&\sum_{j=1}^N\sum_{(k,k')\in {\mathcal N}^\e_j(\Omega)}
\e^{d-1} a_{kk^\prime}{(u_k-u_{k^\prime})^2 }\\
&&+ \sum_{(k,k')\in {\mathcal N}^\e_0(\Omega)} \
\e^{d}a_{kk^\prime} {(u_k-u_{k^\prime})^2 }+\sum_{k\in
Z^\e(\Omega)} \e^d g(k, u_k),
\end{eqnarray}
where
%\begin{equation}
%N^\e_j(\Omega)= \Bigl\{ (k,k')\in (A_j\times A_j)\cap {1\over \e} (\Omega\times\Omega): k-k'\in P^j, k\neq k'\Bigr\},
%\end{equation}
\begin{equation}
{\mathcal N}^\e_j(\Omega)=N_j\cap {1\over \e} (\Omega\times\Omega), j=0,\ldots,N,
\qquad\qquad
Z^\e(\Omega)=\ZZ^d\cap {1\over \e} \Omega.
\end{equation}

The first sum in the energy takes into account all interactions
between points in $A_j$ ({\em hard phases}), which are supposed to
scale differently than those between points in $A_0$ ({\em soft
phase}) or between points in different phases. The latter are contained in the
second sum. The third sum is a zero-order term taking into account
all types of phases with the same scaling .

Note that the first sum may take into account also points in $A_j\setminus C_j$, which form ``islands'' of the hard phase $P^j$-disconnected from the corresponding infinite component. Furthermore, in this energy we may have sites that do not interact at all with hard phases.
%{\bf In order to treat them we introduce a notion of connectedness
%as follows: we say that $k$ and $k'$ are connected if either $k=k'$ or there exists a path $\{k_n\}_{n=0,\ldots,K}$
%such that $k_0=k$, $k_K=k'$ and $(k_n,k_{n-1})\in \bigcup_{j=0}^N N_j$. (where do we use this?)}
%\vfill\eject

\begin{remark}[choice of the parameter space]\rm
The energy is defined on discrete functions parameterized on ${1\over\e}\Omega\cap\ZZ^d$.
The choice of this notation, rather than interpreting $u$ as defined
on $\Omega\cap\e\ZZ^d$ allows a much easier notation for
the coefficients, that in this way are $\e$-independent, rather
than obtained by scaling as in (\ref{ch-c}). 
\end{remark}

\section{Homogenization of perforated discrete domains}
In this section we separately consider the interactions in each infinite connected component of the hard phases introduced above. To that end we fix one of the indices $j$,  with $j>0$, dropping it in the notation of this section (in particular we use the symbol $C$ in place of  $C_j$, etc.), and define the energies
\begin{eqnarray}\label{cale}
{\cal F}_\e(u)={\cal F}_\e\bigl(u,
{1\over\e}\Omega\bigr)=\sum_{(k,k')\in N^\e_C(\Omega)}
\e^{d-1} a_{kk^\prime}{(u_k-u_{k^\prime})^2 }\,,
\end{eqnarray}
where
\begin{equation}
N^\e_C(\Omega)= \Bigl\{ (k,k')\in (C\times C)\cap {1\over \e} (\Omega\times\Omega): k'-k
\in P_k, k\neq k'\Bigr\},
\end{equation}
We also introduce the notation $C^\e(\Omega)= C\cap{1\over\e}\Omega$.

\begin{definition}\label{pci}\rm
The {\em piecewise-constant interpolation} of a function $u:\ZZ^d\cap{1\over\e}\Omega\to\rr^m$, $k\mapsto u_k$ is defined as
$$
u(x)= u_{\lfloor x/\e\rfloor},
$$
where $\lfloor y\rfloor= (\lfloor y_1\rfloor,\ldots, \lfloor y_d\rfloor)$ and $\lfloor s\rfloor$ stands for the integer part of $s$.
The {\em convergence} of a sequence $(u^\e)$ of discrete functions is understood as the $L^1_{\rm loc}(\Omega)$ convergence of these piecewise-constant interpolations. Note that, since we consider local convergence in $\Omega$, the value of $u(x)$ close to the boundary in not involved in the convergence process.
\end{definition}

We prove an extension and compactness lemma with respect to the convergence of piecewise-constant interpolations.

\begin{lemma}[extension and compactness]\label{exc}
Let  $C$ be a $T$-periodic subset of $\ZZ^d$ $P$-connected in the
notation of the previous section, let $u^\e: \ZZ^d\cap{1\over\e}\Omega\to\{+1,-1\}$ be a sequence such that
\begin{equation}
\sup_\e\ \e^{d-1}\#\bigl\{(k,k')\in
N^\e_C(\Omega):u^\e_k\not=u^\e_{k^\prime}\bigr\} <+\infty.
\end{equation}
Then there exists a sequence $\widetilde u^\e:
\ZZ^d\cap{1\over\e}\Omega\to\rr^m$ such that $\widetilde u^\e_k= u^\e_k$ if
$k\in C^\e(\Omega)$ and {\rm dist}$(k, \partial{1\over\e}\Omega)>
c=c(P)$, with $\widetilde u^\e$ converging to some $u\in
BV_{\loc}(\Omega;\{+1,-1\})$ up to subsequences.
\end{lemma}

\Proof 
For a fixed $M\in\NN$ and $j\in \ZZ^d$ we consider the discrete cubes of side length $M$
$$
Q_M(j):=jM+\{0,M-1\}^d.
$$
For each $j$ we also define the cube
$$
Q_{3M}^\prime(j)=\bigcup_{\|i-j\|_\infty\le 1}  Q_M(i),
$$
which is a discrete cube centered in $Q_M(j)$ and with side length $3M$.

For all $\e$ we consider the family 
$$
{\cal Q}_M^\e:=\Bigl\{ Q_M(j): j\in \ZZ^d,\ Q'_{3M}(j)\subset {1\over\e} \Omega\Bigr\}.
$$

We suppose that $M$ is large enough so that if $k, k'\in Q_M(j)\cap C$ then there exists
a $P$-path connecting $k$ and $k'$ contained in $Q'_{3M}(j)$. 
The existence of such $M$ follows easily from the connectedness hypotheses.
Indeed, we may take $M$ as the length of the longest shortest $P$-path connecting 
two points in $C$ with distance not greater than $2\sqrt d$ (in particular belonging
to neighbouring periodicity cubes), and construct such $P$-path by concatenating
a family of those shortest paths.

We define the set of indices
$$
{\cal S}_\e= \bigl\{ j\in \ZZ^d: Q_M(j)\in {\cal Q}_M^\e\hbox{ and } u^\e\hbox{ is not constant on }C\cap Q_M(j)\bigr\}.
$$
By our choice of $M$ if $j\in {\cal S}_\e$ then there exist $k, k'\in Q'_{3M}(j)\cap C$ 
with $k'-k\in P$ such that
$u^\e_k\neq u^\e_{k'}$. Let
$$
K:=\sup_\e\ \e^{d-1}\#\bigl\{(k,k')\in
N^\e_C(\Omega): u^\e_k\not=u^\e_{k^\prime}\bigr\}.
$$
Then we deduce that 
\begin{equation}\label{Kest}
\#{\cal S}_\e\le 3^d K{1\over\e^{d-1}}
\end{equation}
(the factor $3^d$ comes from the fact that $k, k'\in Q'_{3M}(j)$ for $3^d$ possible $j$).

We define $\widetilde u^\e$ as follows
$$
\widetilde u^\e
=\cases{
\hbox{the constant value of }u^\e \hbox{ on } Q_M(j) \cap C &$ \hbox{on } Q_M(j), 
\hbox{ if }Q_M(j)\in {\cal Q}_M^\e
\hbox{ and } j\not\in {\cal S}_\e$\cr
u^\e & elsewhere.}
$$

This will be the required extension. However we will prove the convergence
of $\widetilde u^\e$ as a consequence of the convergence of the functions 
$v^\e$ defined as
$$
v^\e=\cases{\widetilde u^\e&$ \hbox{on } Q_M(j), 
\hbox{ if }Q_M(j)\in {\cal Q}_M^\e
\hbox{ and } j\not\in {\cal S}_\e$\cr
1 & elsewhere.}
$$
By (\ref{Kest}) we have that for fixed $\Omega'\subset\!\subset \Omega$
$$
\|v^\e-\widetilde u^\e\|_{L^1(\Omega')}= O(\e)
$$
(recall that we identify the function with their scaled interpolations in $L^1(\Omega)$).

If the value of $v^\e$ differs on two neighbouring $Q_M(j)$ and $Q_M(j')$ with $\|j-j'\|_1=1$
then, upon taking a suitable larger $M$, we may also suppose that there exist $k,k'\in (Q'_{3M}(j)\cup Q'_{3M}(j))\cap C$ with $k-k'\in P$ and $u^\e_k\not=u^\e_{k^\prime}$.
Arguing as for (\ref{Kest}), we deduce that the number of such $j$ is $O(\e^{1-d})$, so that
$$
{\cal H}^{d-1}(\partial \{v^\e=1\}\cap \Omega')= O(1),
$$
which implies the compactness of the family $(v^\e)$ in $BV_{\rm loc}(\Omega)$.
\qed

\begin{theorem}[homogenization on discrete perforated domains]\label{hdpd}
The energies ${\cal F}_\e$ defined in {\rm(\ref{cale})}
$\Gamma$-converge with respect to the $L^1_{\rm loc} (\Omega)$
topology to the energy
\begin{equation}
{\cal F}_{\hom}(u)=\int_{\Omega\cap\partial^*E} f_{\hom}(\nu)\dHH,
\end{equation}
defined on $u=\chi_E$, $u\in BV(\Omega,\{+1,-1\})$ where the
energy density $f_{\hom}$ satisfies
\begin{eqnarray}\nonumber
f_{\hom}(\xi)&=&\lim_{T\to+\infty} {1\over T^{d-1}} \inf\Bigl\{ 
\sum_{(k,k')\in \widetilde{N}_C(Q^\nu_T)}
 a_{kk^\prime}{(u_k-u_{k^\prime})^2 }: u_k=1 \hbox{ if } k\not\in Q^\nu_T
 \hbox{ and } \langle k,\nu\rangle>0
\\
&&\qquad\qquad u_k=-1 \hbox{ if } k\not\in Q^\nu_T
 \hbox{ and } \langle k,\nu\rangle\le0
\Bigr\},\label{foro}
\end{eqnarray}
where $Q^\nu$ is a cube centered in $0$ and with one side orthogonal to $\nu$, $Q^\nu_T= T\,Q^\nu$, and $\widetilde{N}_C(Q^\nu_T)$ denote all pairs in $(k,k')\in {N}^1_C(\rr^d)$ such that either $k\in Q^\nu_T$ or $k'\in Q^\nu_T$.
\end{theorem}

\Proof In \cite{BP1} this theorem is proved under the
additional assumption that the energies ${\cal F}_\e$ be equi-coercive with respect to the weak $BV$-convergence. This assumption can be substituted by Lemma \ref{exc}. Indeed, if $u^\e$ is a sequence converging to $u$ in $L^1_{\rm loc} (\Omega)$ and with equibounded energies then by Lemma \ref{exc} we may find a sequence $\widetilde u^\e$ coinciding with $u^\e$ on $C^\e(\Omega')$ for every fixed $\Omega'\subset\subset \Omega$ and $\e$ sufficiently small, and converging to some $\widetilde u$ in $BV(\Omega;\{\pm1\})$. Since $\widetilde u^\e= u_\e$ on $C^\e(\Omega')$ we have that $\widetilde u=u$ and ${\cal F}_\e\bigl(\widetilde u^\e, {1\over\e}\Omega'\bigr)={\cal F}_\e\bigl(u^\e, {1\over\e}\Omega'\bigr)$. Then we can give a lower estimate on each $\Omega'$ fixed using the proof of \cite{BP1}, and hence on $\Omega$ by internal approximation. Note that neither the proof of the existence of the limit in (\ref{foro}) therein, nor the construction of the recovery sequences depend on the coerciveness assumption, so that the proof is complete.\qed

\section{Definition of the interaction term}
The homogenization result in Theorem \ref{hdpd} will describe the contribution of the hard phases to
the limiting behavior of energies $F_\e$. We now characterize their interactions with the soft phase.

For all $M$ positive integer and $z_1,\ldots, z_N\in \{+1,-1\}$ we
define the minimum problem
\begin{equation}
\varphi_M(z_1,\ldots, z_N) ={1\over M^{d}} \min\Bigl\{
\sum_{(k,k')\in N_0(Q_M)} a_{k k'}{(v_k-v_{k'})^2 }+\sum_{k\in Z(Q_M)} g(k, v_k): v\in{\cal V}_M\Bigr\},
\end{equation}
where
\begin{equation}
Q_M=\Bigl[-{M\over 2},{M\over 2}\Bigr)^d, \qquad N_0(Q_M)= N_0\cap (Q_M\times Q_M),\qquad Z(Q_M)= \ZZ^d\cap Q_M,
\end{equation}
and the minimum is taken over the set ${\cal V}_M= {\cal V}_M(z_1,\ldots, z_N)$ of all $v$ that are constant on each connected component of $A_j\cap Q_M$ and $v=z_j$
on $C_j$ for $j=1,\ldots N$.

%{\color{red} controllare se/dove si \`e usato la positivit\`a delle interazioni}

\begin{proposition}\label{vf}
There exists the limit $\varphi$ of $\varphi_M$ as $M\to +\infty$.
\end{proposition}

\Proof We first show that 
\begin{equation}
\label{fikm} 
\varphi_{KM} \geq \varphi_{M} \hbox{ for all } K\in\mathbb N. 
\end{equation}
To that end, let $\overline v$ be a minimizer for $\varphi_{KM}(z_1,\ldots, z_N).$ Then we have 
\begin{eqnarray*}
&&\hskip-2cm K^dM^d\varphi_{KM}(z_1,\ldots, z_N)\\
&=& 
\sum_{(k,k')\in N_0(Q_{KM})} a_{k k'}{(\overline v_k-\overline v_{k'})^2 }+\sum_{k\in Z(Q_{KM})} g(k, \overline v_k)\\
&=&\sum_{l\in \ZZ^d\cap Q_K}\Biggl(\sum_{(k,k')\in N_0(Q_M+l M)} a_{k k'}{(\overline v_k-\overline v_{k'})^2 }+\sum_{k\in Z(Q_M+l M)} g(k, \overline v_k)\Biggr)\\
&&+\sum_{(k,k')\in N_0(Q_{KM})\setminus \bigcup_l N_0(Q_M+l M)} 
a_{k k'}{(\overline v_k-\overline v_{k'})^2 }\\
&\geq&\sum_{l\in \ZZ^d\cap Q_K}\Biggl(\sum_{(k,k')\in N_0(Q_M+l M)} a_{k k'}{(\overline v_k-\overline v_{k'})^2 }+\sum_{k\in Z(Q_M+l M)} g(k, \overline v_k)\Biggr). 
\end{eqnarray*} 
Let $\overline l\in\ZZ^d\cap Q_K$ minimize the expression in parenthesis. 
Then we deduce 
\begin{eqnarray*}
K^dM^d\varphi_{KM}(z_1,\ldots, z_N)
&\geq&K^d\Biggl(\sum_{(k,k')\in N_0(Q_M+\overline l M)} a_{k k'}{(\overline v_k-\overline v_{k'})^2 }+\sum_{k\in Z(Q_M+\overline l M)} g(k, \overline v_k)\Biggr), 
\end{eqnarray*} 
from which (\ref{fikm}) follows by taking $v_k=\overline v_{k-\overline l M}$ in the computation of 
$\varphi_{M}(z_1,\ldots, z_N)$.  

We remark that for $L\geq L^\prime$ we have 
\begin{equation}
\label{elle}
L^d\varphi_L\geq (L^\prime)^d \varphi_{L^\prime}-\max|g|\ (L^d-(L^\prime)^d). 
\end{equation}
Hence, fixing $n$, $L$ and $M$, $L\geq M 2^n$, and taking $L^\prime=\lfloor \frac{L}{M 2^n}\rfloor M 2^n$ in (\ref{elle}), we have, using (\ref{fikm}) with $K=\lfloor \frac{L}{M 2^n}\rfloor 2^n$
\begin{eqnarray*}
\varphi_L&\geq& \frac{1}{L^d} \Bigl(\Bigl\lfloor \frac{L}{M 2^n}\Bigr\rfloor M 2^n\Bigr)^d 
\varphi_{\lfloor \frac{L}{M 2^n}\rfloor M 2^n} - \max|g| \Bigl(1-\Bigl(\Bigl\lfloor \frac{L}{M 2^n}\Bigr\rfloor \frac{M 2^n}{L}\Bigr)^d \Bigr)\\ 
&\geq& \Bigl(\Bigl\lfloor \frac{L}{M 2^n}\Bigr\rfloor \frac{M 2^n}{L}\Bigr)^d 
\varphi_M - \max|g| \Bigl(1-\Bigl(\Bigl\lfloor \frac{L}{M 2^n}\Bigr\rfloor \frac{M 2^n}{L}\Bigr)^d \Bigr).
\end{eqnarray*}  
Letting $L\to+\infty$ we then obtain 
$$\liminf_{L\to+\infty} \varphi_L\geq \varphi_M$$ 
and the thesis by taking the upper limit in $M$. 
\qed

Let $R$ be defined by 
\begin{equation}
\label{erre} 
R=\max\{ |k-k^\prime| : \ k, k' \in A_j\setminus C_j \ P^j\hbox{-connected}, \ j=1,\ldots, N\}
\end{equation}
and for all $M$ positive integer we set 
\begin{equation}
\label{diemme} 
D_M=\bigcup_{j=1}^N 
\bigcup \bigl\{B: \ B \hbox{ a } P^j\hbox{-connected components of } A_j\setminus C_j \hbox{ not intersecting } Q_{M-R}\bigr\}.
\end{equation}
For all $z_1,\ldots, z_N\in \{+1,-1\}$ we define 
\begin{eqnarray*}
\widetilde\varphi_M(z_1,\ldots, z_N)&=&
{1\over M^{d}} \min\Bigl\{
\sum_{(k,k')\in N_0(Q_M)} a_{k k'}{(v_k-v_{k'})^2 }+\sum_{k\in Z(Q_M)} g(k, v_k): \\
&&\qquad\qquad\qquad v\in{\cal V}_M, v_k=1 \hbox{ if } k\in D_M\Bigr\}.
\end{eqnarray*}

\begin{proposition}\label{fikm2}
There is a positive constant $c$ independent of $M$ such that 
\begin{equation}
\label{tilde} 
\widetilde\varphi_M\geq \varphi_M\geq \widetilde\varphi_M-\frac{c}{M}.
\end{equation}
\end{proposition}
\Proof
The first inequality is trivial. To prove the second, let $\overline v$ be a minimizer for $\varphi_M(z_1,\ldots, z_N)$ and define $v$ by 
$$v_k=\cases{
1 & if  $k\in D_M$ \cr
\overline v_k & otherwise.} 
$$
Using $v$ as a test function for $\widetilde\varphi_M(z_1,\ldots, z_N)$, we obtain 
\begin{eqnarray*}
M^d\widetilde\varphi_M(z_1,\ldots, z_N)&\leq& 
\sum_{(k,k')\in N_0(Q_M), \ k,k'\not\in D_M} a_{k k'}{(v_k-v_{k'})^2 }+\sum_{k\in Z(Q_M)\setminus D_M} g(k, v_k)\\
&&+2\sum_{(k,k')\in N_0(Q_M), \ k\in D_M} a_{k k'}{(v_k-v_{k'})^2 }+\sum_{k\in Z(Q_M)\cap D_M} g(k, v_k)\\
&\leq&
\sum_{(k,k')\in N_0(Q_M), \ k,k'\not\in D_M} a_{k k'}{(\overline v_k-\overline v_{k'})^2 }+\sum_{k\in Z(Q_M)\setminus D_M} g(k, \overline v_k)\\
&&+\sum_{(k,k')\in N_0(Q_M), \ k\in D_M} a_{k k'}+\sum_{k\in Z(Q_M)\cap D_M} g(k, 1)\\
&\leq&
M^d \varphi_M(z_1,\ldots, z_N)
+\# D_M \ \# P_0 \ \max a_{ij} + \# D_M \ 2 \max |g|.  
\end{eqnarray*}
Since $\# D_M\leq 2^d M^{d-1} R,$ the thesis follows with 
$c=2^d R( \# P_0 \ \max a_{ij} + \ 2 \max |g|)$. 
\qed

%
%
%\begin{remark}\label{frame}\rm
%Let $u^M\in {\cal V}_M$ be a sequence such that
%$$
%\lim_M{1\over M^d}\Bigl(\sum_{(k,k')\in N_0(Q_M)} f(k,k', {u^M_k-u^M_{k'}})+\sum_{k\in Z(Q_M)} g(k, u^M_k)\Bigr)=\varphi(z_1,\ldots, z_N) $$
%then for every sequence of constants $R_M= o(M)$ we have
%$$
%\lim_M{1\over M^d}
%\sum_{\buildrel{\scriptstyle k,k'\in Q_M\setminus Q_{M-R_M}}\over{\scriptstyle k-k'\in P_0}}  |u^M_k-u^M_{k'}|^p=0.
%$$
%Indeed, otherwise taking $u^M$ as test function for the problem defining $\varphi_{M-R_M}(z_1,\ldots, z_N)$,
%we would obtain
%$$
%\limsup_M\varphi_{M-R_M}(z_1,\ldots, z_N)<\varphi(z_1,\ldots, z_N),
%$$
%which is a contradiction.
%\end{remark}
%

\section{Statement of the convergence result}
We now have all the ingredients to characterize the asymptotic behavior of $F_\e$ defined in (\ref{fep}).

\begin{definition}[multi-phase discrete-to-continuum convergence]
We define the {\em convergence}
\begin{equation}\label{conva}
u^\e\to (u^1,\ldots, u^N)
\end{equation}
as the $L^1_{\rm loc} (\Omega;\rr^m)$ convergence $\widetilde u^\e_j\to u^j$ of the extensions
of the restrictions of $u^\e$ to $C_j$ as in Lemma {\rm\ref{exc}}, which is a compact convergence as ensured by that lemma.
\end{definition}

The total contribution of the hard phases will be given separately by the contribution on the infinite connected components and the finite ones.
The first one is obtained by
computing independently the limit relative to the energy restricted to each component
\begin{eqnarray}\label{calej}
\mathcal F^j_\e(u)=\sum_{(k,k')\in N^\e_j(\Omega)} \e^{d-1} 
a_{k k'}{(v_k-v_{k'})^2 }\,,
\end{eqnarray}
where
\begin{equation}N^\e_j(\Omega)= N^\e_{C_j}(\Omega)=
\Bigl\{ (k,k')\in (C_j\times C_j)\cap {1\over \e} (\Omega\times\Omega): k-k'\in P_k^j, k\neq k'\Bigr\},
\end{equation}
which is characterized by Theorem \ref{hdpd} as
\begin{equation}\label{fjh}
\mathcal F^j_{\hom} (u)=\int_{\Omega\cap\partial^*\{u=1\}} f^j_{\hom}(\nu)\dHH.
\end{equation}

%\begin{equation}
%A_j\setminus C_j=\bigcup_{l\in I_j}(A^j_l+ T\ZZ^d),
%\end{equation}
%where, due to the periodicity of the media, $l$ runs over a finite set of indices
%$I_j$, and $A^j_l+ T\ZZ^d$ and $A^j_{l'}+ T\ZZ^d$ are $P^j$-disconnected if $l\neq l'$.
%To each such $A^j_l$ we associate the minimum value
%\begin{equation}\label{mji}
%m^j_l=\min\Bigl\{\sum_{k,k'\in A^j_l, k-k'\in P^j} f(k,k', {z_k-z_{k'}}): z:A^j_l\to \rr^m\Bigr\}.
%\end{equation}
%Note that we have no boundary conditions for the test functions $z$. The total contribution of the disconnected components will simply give the additive constant $m|\Omega|$, where
%\begin{equation}\label{m}
%m={1\over T^d} \sum_{j=1}^N\sum_{l\in I_j} m^j_l.
%\end{equation}

In the previous section we have introduced the energy density $\varphi$, which describes the interactions between the hard phases. Taking all contribution into account, we may state the following convergence result.

\begin{theorem}[double-porosity homogenization]\label{dph}
Let $\Omega$ be a Lipschitz bounded open set, and let $F_\e$ be defined by {\rm(\ref{fep})} with the notation of Section {\rm\ref{Not}}. Then there exists the $\Gamma$-limit of $F_\e$ with respect to the convergence {\rm(\ref{conva})}
and it equals
\begin{equation}\label{maineq}
F_{\hom}(u^1,\ldots, u^N)=\sum_{j=1}^N\int_{\Omega\cap\partial^*\{u^j=1\}} f^j_{\hom}(\nu)\dHH+
\int_\Omega \varphi(u^1,\ldots, u^N)\dx
\end{equation}
on functions $u=(u^1,\ldots, u^N)\in (BV(\Omega;\{1,-1\}))^N$,
where $\varphi$ is defined in Proposition {\rm \ref{vf}}, $f^j_{\hom}$ are defined by  {\rm(\ref{fjh})}.
\end{theorem}

Note that there is no contribution of the finite connected components of $A_j$.

The proof of this result will be subdivided into a lower and an upper bound. 

%{\color{red} riscrivere $\varphi$ come somma di misure di insiemi?}

\begin{remark}[non-homogeneous lower-order term]\rm
In our hypotheses the lower-order term $g$ depends on the fast variable $k$, which is integrated out in the limit.
We may easily include a measurable dependence on the slow variable $\e k$, by assuming $g= g(x,k,z)$ a Carath\'eodory function (this covers in particular the case $g=g(x,z)$) and substitute the last sum in (\ref{fep}) by
$$
\sum_{k\in Z^\e(\Omega)} \e^d g(\e k,k, u_k).
$$
Correspondingly, in Theorem \ref{dph} the integrand in the last term in (\ref{maineq}) must be substituted by
$\varphi(x,u^1,\ldots, u^N)$, where the definition of this last function is the same but taking $g(x,k,z)$ in place of $g(k,z)$, so that $x$ simply acts as a parameter.
\end{remark}

\subsection{Proof of the lower bound} 
Let $u^\e\to(u^1,\ldots, u^N)$ be such that $F_\e(u^\e)\le c<+\infty$. 
Fixed $M\in\NN$, 
we introduce the notation 
$$J^\e_M= \{z\in\ZZ^d : \ Q_M+zM\subset \frac{1}{\e}\Omega\},$$
$$R^\e=\mathcal N^\e_0(\Omega)\setminus\bigcup_{z\in J^\e_M}\mathcal N^\e_0(Q_M+zM),
$$
$$
S^\e= Z^\e(\Omega) \setminus\bigcup_{z\in J^\e_M} Z(Q_M+zM),
$$
and write 
$$F_\e(u^\e)=\sum_{j=1}^N I^\e_j+I\!I^\e+I\!I\!I^\e+IV^\e+V^\e,
$$ 
where 
$$I^\e_j= \mathcal F^j_\e(u),$$
$$I\!I^\e= \sum_{j=1}^N\sum_{(k,k')\in \mathcal N^\e_j(\Omega)\setminus(C_j\times C_j)} \e^{d-1} 
a_{k k'}{(v_k-v_{k'})^2 },$$
$$I\!I\!I^\e= \sum_{z\in J^\e_M}\e^d \Biggl( \sum_{(k,k')\in \mathcal N^\e_0(Q_M+zM)} a_{k k'}{(v_k-v_{k'})^2 }+\sum_{k\in Z(Q_M+zM)} g(k, v_k)\Biggr),$$
$$IV^\e=\sum_{(k,k')\in R^\e} \e^d a_{k k'}{(v_k-v_{k'})^2 },
$$
and
$$V^\e=\sum_{k\in S^\e} \e^d g(k, v_k).
$$

Note that  
\begin{equation}
\label{unno} 
I\!I^\e\geq 0,\qquad IV^\e
%\geq 0
\geq -c/M+o(1), \qquad V^\e \geq - \max|g| \Bigl(\Bigl|\Omega\setminus \e^d\bigcup_{z\in J^\e_M}(Q_M+zM)\Bigr|+o(1)\Bigr),
\end{equation}
where we have taken into account that the interactions in $IV^\e$ may be negative.
and 
\begin{equation}
\label{ddue}
\liminf_{\e\to 0} \sum_{j=1}^N I^\e_j \geq  \sum_{j=1}^N \liminf_{\e\to 0} I^\e_j\geq 
\sum_{j=1}^N\int_{\Omega\cap\partial^*\{u^j=1\}} f^j_{\hom}(\nu)\dHH.
\end{equation}

It remains to estimate $I\!I\!I^\e$. To that end, we introduce the set of indices 
$$\Lambda_M^\e=\{z\in J_M^\e : \ u^\e \hbox{ constant on every connected component of } 
A_j\cap (Q_{3M}+zM), j=1,\ldots, N\}.$$
Note that 
\begin{equation}
\label{lambda}
\#(J_M^\e\setminus \Lambda_M^\e)\leq \frac{c_M}{\e^{d-1}}.
\end{equation}
We then write 
\begin{eqnarray*}
I\!I\!I^\e&=&\sum_{z\in \Lambda^\e_M}\e^d \Biggl( \sum_{(k,k')\in \mathcal N^\e_0(Q_M+zM)} a_{k k'}{(v_k-v_{k'})^2 }+\sum_{k\in Z(Q_M+zM)} g(k, v_k)\Biggr)\\
&&+\sum_{z\in J^\e_M\setminus \Lambda_M^\e}\e^d \Biggl( \sum_{(k,k')\in \mathcal N^\e_0(Q_M+zM)} a_{k k'}{(v_k-v_{k'})^2 }+\sum_{k\in Z(Q_M+zM)} g(k, v_k)\Biggr)\\
&\geq& \sum_{z\in \Lambda^\e_M}\e^d M^d \varphi_M(u^\e_1,\ldots, u_N^\e)- c
\e^d  M^d \max(|g|
%{\color{red}
+|a_{kk'}|) \# (J^\e_M\setminus \Lambda_M^\e),   
\end{eqnarray*}
where $u_j^\e$ is the constant value taken by $u^\e$ on $(Q_M+zM)\cap C_j$. Here we suppose $M$ large enough so that the connected component of $C_j$ containing $(Q_M+zM)\cap C_j$ is connected in $Q_{3M}+zM$. 
We set 
$$U^\e=\sum_{z\in \Lambda_M^\e} (u_1^\e,\ldots, u_N^\e)\chi_{Q_M+zM}$$
and $\varphi_M(0,\ldots, 0)=0$. 
Note that $U^\e\to U:=(u^1,\ldots, u^N)$ in $L^1(\Omega)^N$, so that 
\begin{equation}
\label{ttre}
\liminf_{\e\to 0} I\!I\!I^\e\geq \liminf_{\e\to 0}\Bigl(\int_{\Omega} \varphi_M(U^\e)\, dx-\e \max
%{\color{red}
|g| c_M M^d\Bigr)=\int_{\Omega} \varphi_M(U)\, dx
\end{equation} 
by the Lebesgue dominated convergence theorem and the estimate (\ref{lambda}).  

Summing up the inequalities (\ref{unno}), (\ref{ddue}) and (\ref{ttre}), we get 
\begin{equation}
\label{liminf}
\liminf_{\e\to 0} F_\e(u^\e)\geq \sum_{j=1}^N\int_{\Omega\cap\partial^*\{u^j=1\}} f^j_{\hom}(\nu)\dHH+\int_{\Omega} \varphi_M(U)\, dx.
\end{equation}
The lower bound inequality then follows by taking the limit as $M\to+\infty$, using Proposition \ref{vf} and the Lebesgue dominated convergence theorem. 

\subsection{Proof of the upper bound} 
We fix $U=(u^1,\ldots, u^N)\in BV(\Omega;\{1,-1\})^N$. For every $j=1,\ldots,N$ we choose 
$u^{j,\e}\to u^j$ a recovery sequence for $\mathcal F^j_{\rm hom}(u^j)$. We tacitly extend all functions defined on $Z^\e(\Omega)$ to the whole $\ZZ^d$ with the value $+1$ outside 
$Z^\e(\Omega)$. This does not affect the value of the energies, but allows to rigorously define some sets of indices $z$ in the sequel.

We fix $M\in \NN$ large enough. Similarly as in the previous section, we introduce the sets of indices 
$$
\widetilde J_M^\e=\Bigl\{ z\in\ZZ^d: (Q_M+zM)\cap {1\over\e}\Omega\neq\emptyset\Bigr\},
$$
$$\Lambda_M^{j,\e}=\{z\in J_M^\e : \ u^\e \hbox{ constant on every connected component of } 
A_j\cap (Q_{3M}+zM)\},$$ 
% qui sarebbe meglio definire J tilde epsilon che tiene conto di tutti i quadrati che intersecano Omega
and remark the estimate 
\begin{equation}
\label{lambdaj}
\sum_{j=1}^N\#(\widetilde J_M^\e\setminus \Lambda_M^{j,\e})\leq \frac{c_M}{\e^{d-1}}.
\end{equation}
Note that if $z\in \bigcap_{j=1}^N\Lambda_M^{j,\e}$ then $u^{j,\e}=:u^{j,\e,z}$ is constant on $C_j\cap(Q_{M}+zM)$ for $j=1,\ldots, N$. 
Let $v^{\e,z}$ be a minimizer for $\widetilde\varphi_M(u^{1,\e,z},\ldots,u^{N,\e,z})$. 

We define 
$$u^\e_k= \cases{
u^{j,\e}_k& if $k\in C_j, \ j=1,\ldots, N$\cr
v^{\e,z}(k-zM)& if $k\in Q_M+zM$ and $z\in \bigcap_{j=1}^N\Lambda_M^{j,\e}$\cr 
1& otherwise. 
}
 $$

We first estimate the energy on the strong connections. By the definition of 
$u^{j,\e}$ we have, for all $j=1,\ldots, N$
%stima su N_j 
\begin{equation}
\lim_{\e\to 0}\sum_{(k,k')\in \mathcal N^\e_j(\Omega)\cap(C_j\times C_j)} \e^{d-1} a_{k k'}{(u^\e_k-u^\e_{k'})^2 }= 
\mathcal F^j_{\rm hom}(u^j), 
\end{equation}
since $u^{\e}=u^{j,\e}$ on $C_j$. On the strong connections between points not in the infinite connected components $C_j$ we have 
\begin{equation}
\sum_{(k,k')\in \mathcal N^\e_j(\Omega)\setminus(C_j\times C_j)} \e^{d-1} a_{k k'}{(u^\e_k-u^\e_{k'})^2 }= 0
\end{equation}
since $u^\e$ is constant on every connected component of $A_j\setminus C_j$. Note that here we have used the condition that $v^{\e,z}=1$ on $D_M$ in the definition of $\widetilde\varphi_M$. 

%stima per phi
We then examine the contribution due to the interaction between weak connections and the term $g$. We first look at the contributions on the cubes in the sets $\Lambda_M^{j,\e}$,
where we can use the definition of $\widetilde \varphi_M$:
for every $z\in \bigcap_{j=1}^N\Lambda_M^{j,\e}$ we have 
\begin{equation}\label{vez} 
%\sum_{z\in \bigcap_{j=1}^N\Lambda_M^{j,\e}}\e^d\Bigl( 
\sum_{(k,k')\in \mathcal N^\e_0(Q_M+zM)} a_{k k'}{(u^\e_k-u^\e_{k'})^2 }+\sum_{k\in Z(Q_M+zM)} g(k, u^\e_k)=\widetilde\varphi_M(u^{1,\e,z},\ldots,u^{N,\e,z}).
%\Bigr)
\end{equation}
The contributions interior to all other cubes in $ \widetilde J_M^\e$ sums up to
\begin{eqnarray*}
&&\sum_{z\not\in \bigcap_{j=1}^N\Lambda_M^{j,\e}}\e^d\Bigl( 
\sum_{(k,k')\in \mathcal N^\e_0(Q_M+zM)} a_{k k'}{(u^\e_k-u^\e_{k'})^2 }+\sum_{k\in Z(Q_M+zM)} g(k, u^\e_k)
\Bigr)\\
&\leq&\e^d M^d (\# P_0 \max a_{il}+\max|g|)
\sum_{j=1}^N\#(J_M^\e\setminus \Lambda_M^{j,\e})\\
&\leq&\e M^d c'_M +o(1)
\end{eqnarray*}
by (\ref{lambdaj}) and the fact that the boundary of $\Omega$ has zero measure.
Finally, the contribution due to the weak connection across the boundary of neighbouring
cubes is given by
\begin{eqnarray*}
&&\sum_{z\neq z' \in \bigcap_{j=1}^N\Lambda_M^{j,\e}}\e^d
\sum_{(k,k')\in \mathcal N^\e_0(\Omega),
k\in Q_M+zM, k'\in Q_M+z'M} a_{k k'}{(u^\e_k-u^\e_{k'})^2 }\\
&\leq&\e^d M^{d-1} \#J_M^\e \# P_0 \max a_{il}
 \leq  \# P_0 \max a_{il} \frac{|\Omega|}{M}. 
\end{eqnarray*}

From the inequalities above, we obtain 
$$ 
\limsup_{\e\to 0} F_\e(u^\e)\leq \sum_{j=1}^N \mathcal F_{\rm hom}^j(u^j) +
\int_{\Omega} \widetilde\varphi_M(u^1,\ldots,u^N)\, dx+\# P_0 \max a_{il}\frac{|\Omega|}{M}
$$
The thesis is then obtained by letting $M\to+\infty$ and using Propositions \ref{fikm2} and \ref{vf}.
\qed

\section{Examples}

In the pictures in the following examples weak connections are denoted by a dashed line,
strong connections by a continuous line.

\subsection{One-dimensional examples}
In this section we consider very easy one-dimensional examples, highlighting the 
possibility of a double-porosity behaviour if long-range interactions are allowed,
contrary to the continuum case.
We use a slightly different notation that that followed hitherto, with the sums depending only on one index.
The factor $1/4$ is just a normalization factor since $(u_i-u_j)^2$ is always a multiple of $4$.

\begin{figure}[htbp]
\begin{center}
\includegraphics[width=.6\textwidth]{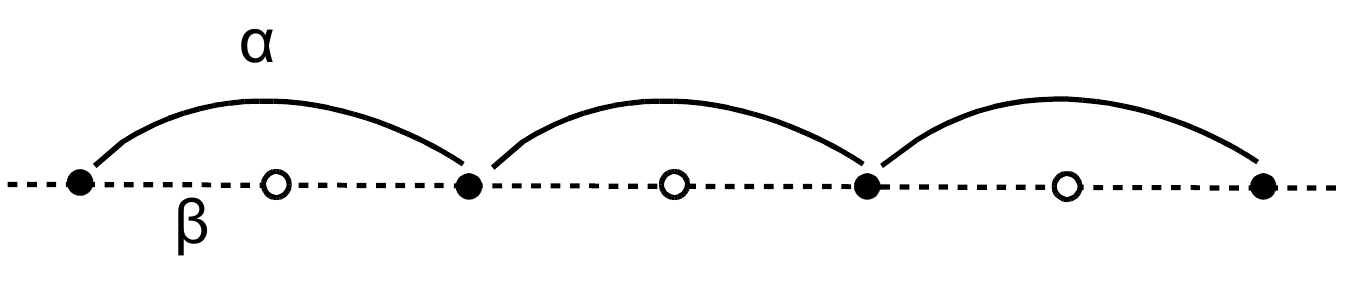}
\caption{weak inclusions in one dimension}
\label{Example1}
\end{center}
\end{figure}
\begin{example}\label{examp-1}\rm
We consider a system of weak nearest-neighbour interactions and strong next-to-nearest neighbour interactions
on the odd lattice (see Fig.~\ref{Example1}); namely, 
$$
F_\e(u)={\beta\over 4}\sum_{i=1}^{N_\e}\e(u_i-u_{i-1})^2+{\alpha\over 4}\sum_{j=1}^{N_\e/2-1}(u_{2j+1}-u_{2j-1})^2+\sum_{i=1}^{N_\e}\e g(u_i)
$$
where we assume that $\Omega=[0,1]$, 
$N_\e=1/\e\in 2\NN$. In this case $N=1$, $A_1=C_1=1+2\NN$, $A_0=2\NN$.

The $\Gamma$-limit is
\begin{eqnarray*}
F_{\rm hom}(u)&=&\alpha\,\# S(u)+ {1\over 2}\int_0^1g(u)\,dx+{1\over 2}\int_0^1\min\{g(u),g(-u)+2\beta\}\,dx\\
&=&\alpha\,\# S(u)+ \int_0^1g(u)\,dx-{1\over 2}\int_0^1\max\{0,g(u)-g(-u)-2\beta\}\,dx
\end{eqnarray*}

The last term favours states with the same value on $A_0$ and $A_1$ if the integrand is $0$ 
and of opposite sign if the integrand is positive. Note that this is always the case if we have a strong enough
`antiferromagnetic' nearest-neighbour interaction; i.e., $\beta$ is negative and $2|\beta|>|g(1)-g(-1)|$.
\end{example}

\begin{figure}[htbp]
\begin{center}
\includegraphics[width=.6\textwidth]{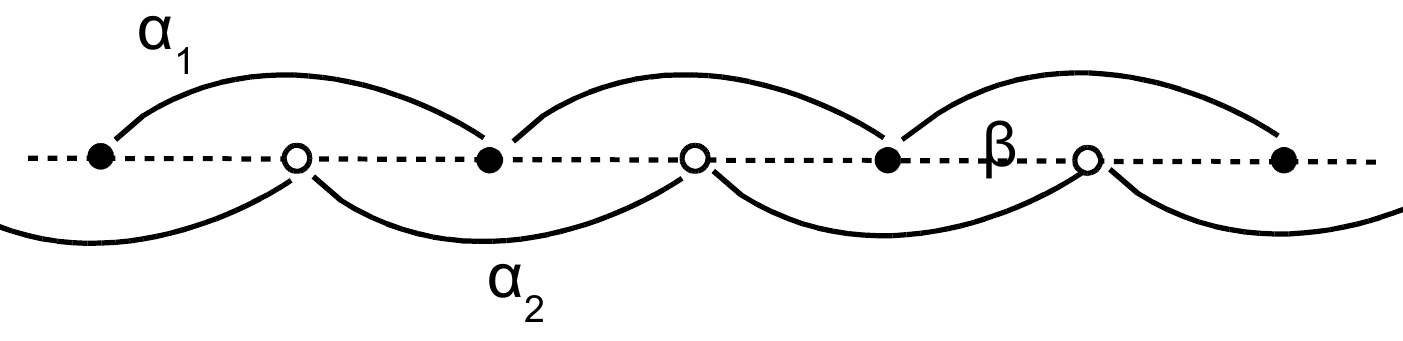}
\caption{interacting sublattices in one dimension}
\label{Example2}
\end{center}
\end{figure}

\begin{example}\label{examp-2}\rm
We consider a system of weak nearest-neighbour interactions and strong next-to-nearest neighbour interactions (see Fig.~\ref{Example2}); namely, 
$$
F_\e(u)={\beta\over 4}\sum_{i=1}^{N_\e}\e(u_i-u_{i-1})^2+{\alpha_1\over 4}\sum_{j=1}^{N_\e/2-1}(u_{2j+1}-u_{2j-1})^2+{\alpha_2\over 4}\sum_{j=0}^{N_\e/2-1}(u_{2j+2}-u_{2j})^2+\sum_{i=1}^{N_\e}\e g(u_i)
$$
where we assume that $N_\e=1/\e\in 2\NN$. In this case $N=2$, $A_1=C_1=1+2\NN$, $A_2=C_2=2\NN$, $A_0=\emptyset$.

The $\Gamma$-limit is
\begin{eqnarray*}
F_{\rm hom}(u^1,u^2)=\alpha_1\,\# S(u^1)+\alpha_2\,\# S(u^2)+ {1\over 2}\int_0^1g(u^1)\,dx+{1\over 2}\int_0^1g(u^2)\,dx+{\beta\over 4}\int_0^1(u^2-u^1)^2.
\end{eqnarray*}
Note that, since $A_0=\emptyset$ we have no optimization in the interacting term, which then is just the pointwise
limit of the nearest-neighbour interactions.
Note moreover that in the case $\beta=0$ the interactions are completely decoupled.
\end{example}

 \begin{figure}[htbp]
\begin{center}
\includegraphics[width=.6\textwidth]{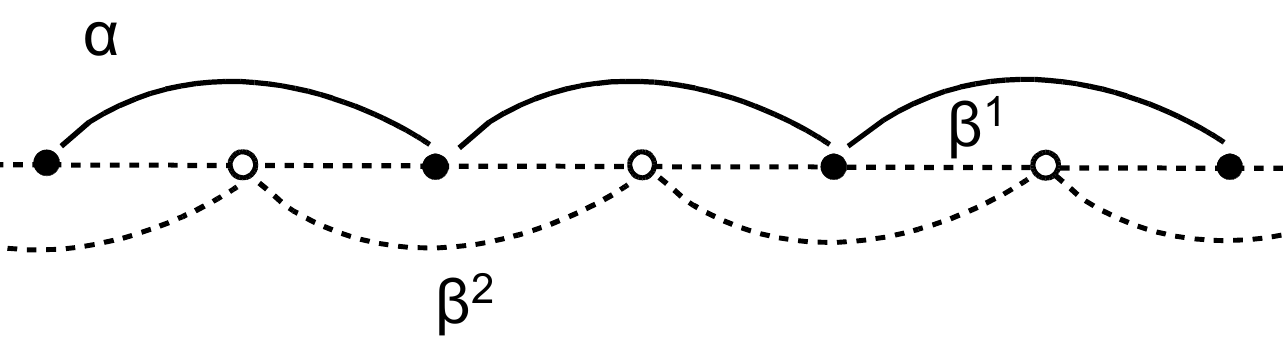}
\caption{interacting weak and strong sublattices in one dimension}
\label{Example2bis}
\end{center}
\end{figure}

\begin{example}\label{examp-3}\rm
We consider the same pattern of interactions as in the previous example, but
with only strong connections on the odd lattice as in Example \ref{examp-1} (see Fig.~\ref{Example2bis}); i.e.,
with 
$$
F_\e(u)={\beta^1\over 4}\sum_{i=1}^{N_\e}\e(u_i-u_{i-1})^2
+{\beta^2\over 4}\sum_{j=0}^{N_\e/2-1}\e(u_{2j+2}-u_{2j})^2
+{\alpha\over 4}\sum_{j=1}^{N_\e/2-1}(u_{2j+1}-u_{2j-1})^2+\sum_{i=1}^{N_\e}\e g(u_i).
$$
In this case we have three possibilities:

1) the minimizing values on the even lattice agree with those on the odd lattice 
(ferromagnetic overall behaviour),

2) the  minimizing values on the even lattice disagree with those on the odd lattice 
(antiferromagnetic overall behaviour),

3) the values on the even lattice alternate (antiferromagnetic behaviour on the weak lattice).

 The value of $\varphi$ is obtained by optimizing on these three possibilities; i.e.,
$$
\varphi(u)=\min\Bigl\{g(u), {g(u)+g(-u)\over 2}+\beta^1, {3g(u)+g(-u)\over 4}+{\beta^1+\beta^2\over 2}\Bigl\},
$$
and we have
$$
F_{\rm hom}(u)=\alpha\,\# S(u)+\int_0^1\varphi(u)\dx.
$$
\end{example}

\begin{figure}[htbp]
\begin{center}
\includegraphics[width=.6\textwidth]{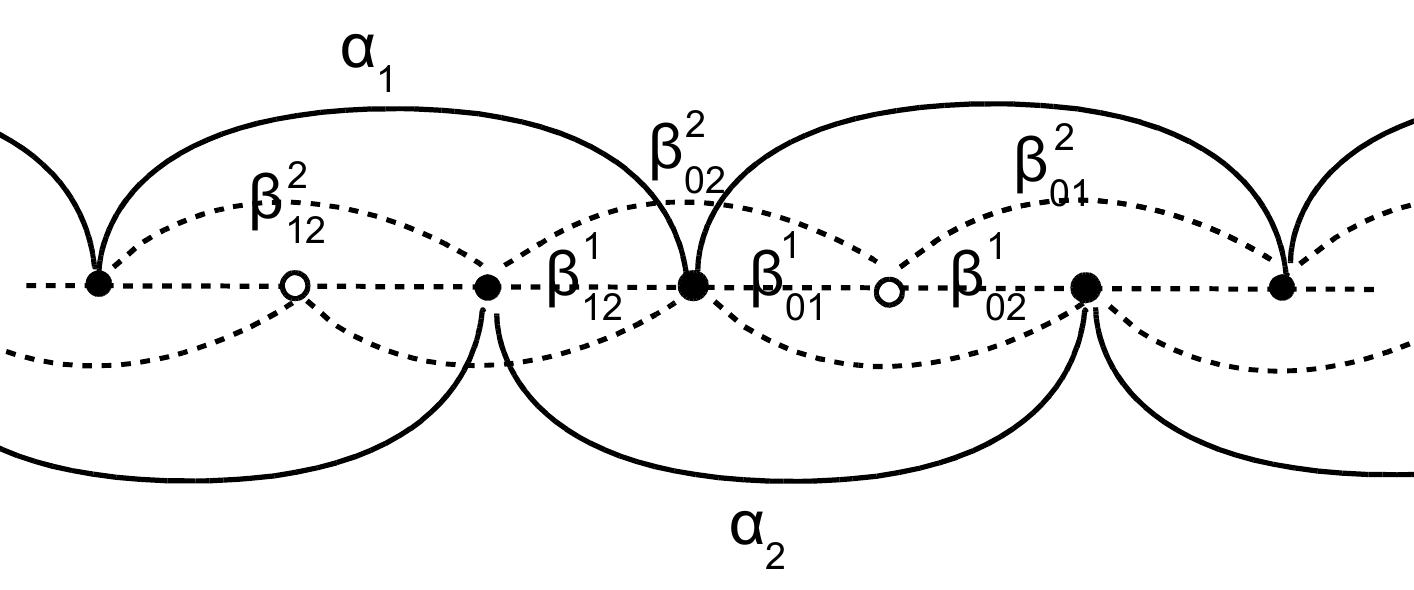}
\caption{Third-neighbour hard phases}
\label{Example3}
\end{center}
\end{figure}

\begin{example}\rm
In the system described in Fig.~\ref{Example3} involving strong third-neighbour interactions, we have two strong components, and a $\Gamma$-limit obtained by minimization of the nearest and next-to nearest neighbours.
Using the same notation of the previous examples for the coefficients as in Fig.~\ref{Example3}, we can write the limit as
\begin{eqnarray*}
F_{\rm hom}(u^1,u^2)=\alpha_1\,\# S(u^1)+\alpha_2\,\# S(u^2)+\int_0^1\varphi(u^1,u^2)\dx,
\end{eqnarray*}
and
\begin{eqnarray*}
\varphi(u^1,u^2)&=&{1\over 3}(g(u^1)+g(u^2))+{1\over 4}\beta^2_{12} (u^2-u^1)^2\\
&&+{1\over 3}\min\Bigl\{{1\over 4}\Bigl((\beta^1_{01}+\beta^2_{01})(v-u^1)^2+(\beta^1_{02}+\beta^2_{02})(v-u^2)^2\Bigr)
+g(v): v\in\{-1,1\}\Bigr\}.
\end{eqnarray*}
\end{example}

\subsection{Higher-dimensional examples}
In the following examples we go back to the notation used in the statement of the main result.
The normalization factor $1/8$ takes into account that each pair of nearest neighbours is accounted 
for twice.

\begin{figure}[htbp]
\begin{center}
\includegraphics[width=.25\textwidth]{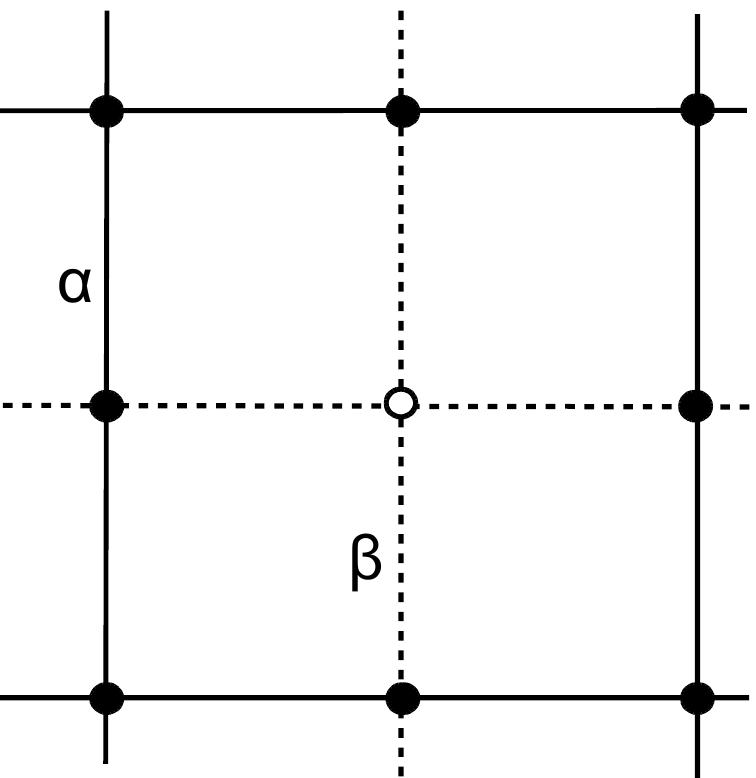}
\caption{a nearest-neighbour system with soft inclusions}
\label{Example5}
\end{center}
\end{figure}
\begin{example}\rm
We consider a nearest-neighbour system in two-dimension in which 
$A_0=2\ZZ^2$ and strong/weak interactions of the form
$$
{1\over 8}\alpha (u_k-u_{k'})^2,\qquad {1\over 8}\e\beta (u_k-u_{k'})^2,
$$
respectively (see Fig.~\ref{Example5}). 
In this case we have
$$
F_{\rm hom}(u)=
{1\over 2}\alpha\int_{S(u)\cap\Omega}\|\nu_u\|_1d\HH^1+
\int_\Omega \varphi(u)\dx,
$$
where
$$
\varphi(u)=\min\Bigl\{g(u), {3g(u)+g(-u)\over 4}+\beta\Bigr\}.
$$
\end{example}

\begin{figure}[htbp]
\begin{center}
\includegraphics[width=.25\textwidth]{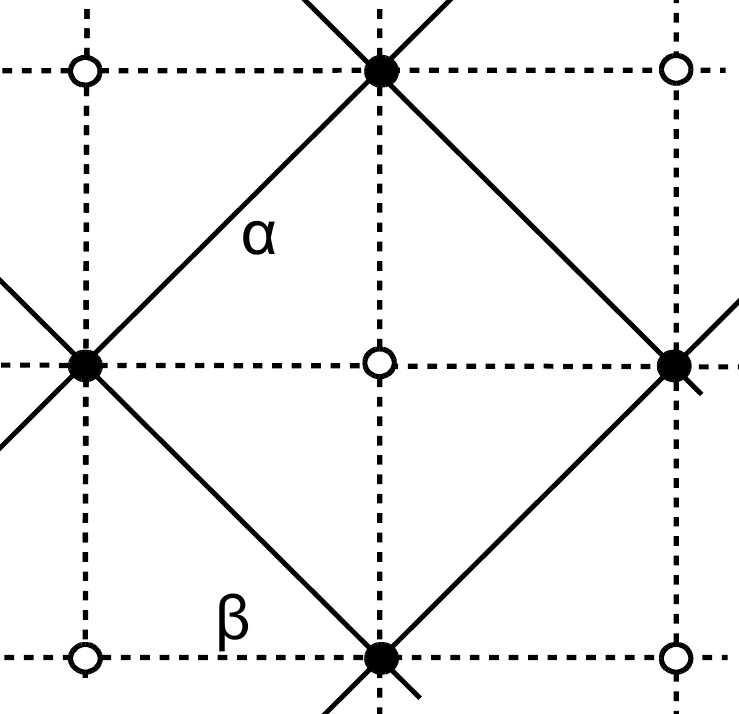}
\caption{a lattice with weak nearest-neighbour interactions}
\label{Example4}
\end{center}
\end{figure}

\begin{example}\rm
In this example we consider strong interactions on a lattice of next-to-nearest neighbours 
as in Fig.~\ref{Example4}, and weak nearest-neighbour interactions on the square lattice, of the form
$$
{1\over 8}\alpha (u_k-u_{k'})^2,\qquad {1\over 8}\e\beta (u_k-u_{k'})^2
$$
respectively (the factor $8$ taking into account that each pair is accounted for twice). 
We only have one strong component, and with this choice of coefficients we have
$$
F_{\rm hom}(u) =\alpha\int_{\Omega\cap\partial\{u=1\}} \|\nu\|_\infty d\HH^1+\int_\Omega \varphi(u)\dx,
$$
where
$$\varphi(u)=\min \Bigl\{g(u),{1\over 2}(g(u)+g(-u))+\beta\Bigr\}.
$$
\end{example}

\begin{figure}[htbp]
\begin{center}
\includegraphics[width=.25\textwidth]{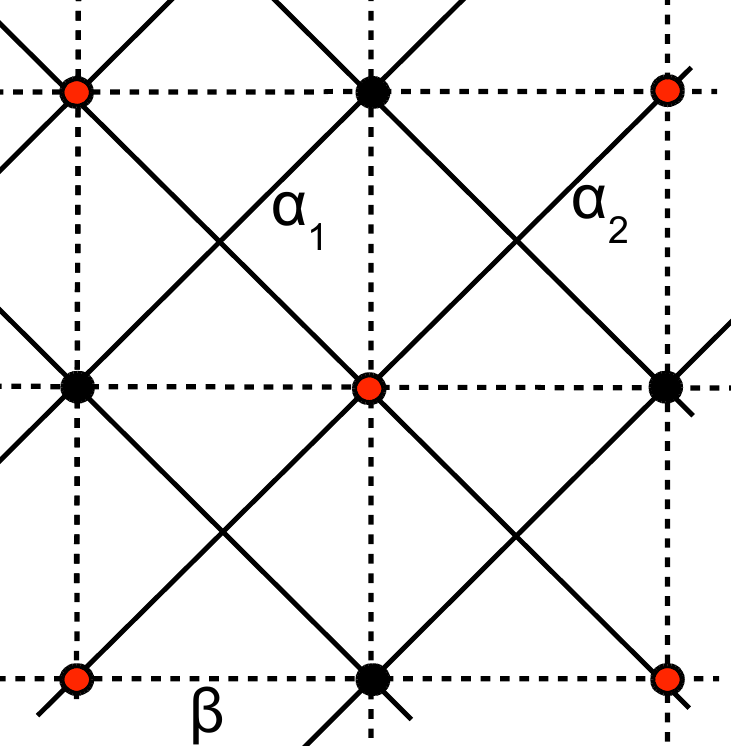}\qquad \qquad\qquad \includegraphics[width=.25\textwidth]{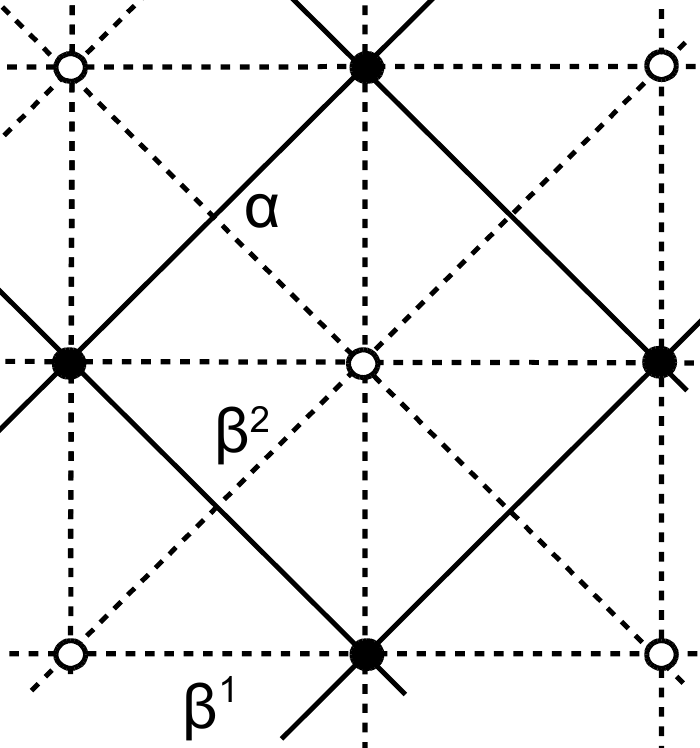}
\caption{two-dimensional interacting sublattices}
\label{Example4bis}
\end{center}
\end{figure}

\begin{example}\rm
We include just the pictorial description of two more two-dimensional systems with a limit with two parameters
(the first one in Fig.~\ref{Example4bis}), 
and with one limit parameter but with
the possibilities of an oscillating behaviour 
on the weak lattice (the second one in Fig.~\ref{Example4bis}),
analogous to the one-dimensional Example \ref{examp-2} and  Example \ref{examp-3},
respectively.
\end{example}

\begin{figure}[h!]
\centerline{\includegraphics [width=2.5in]{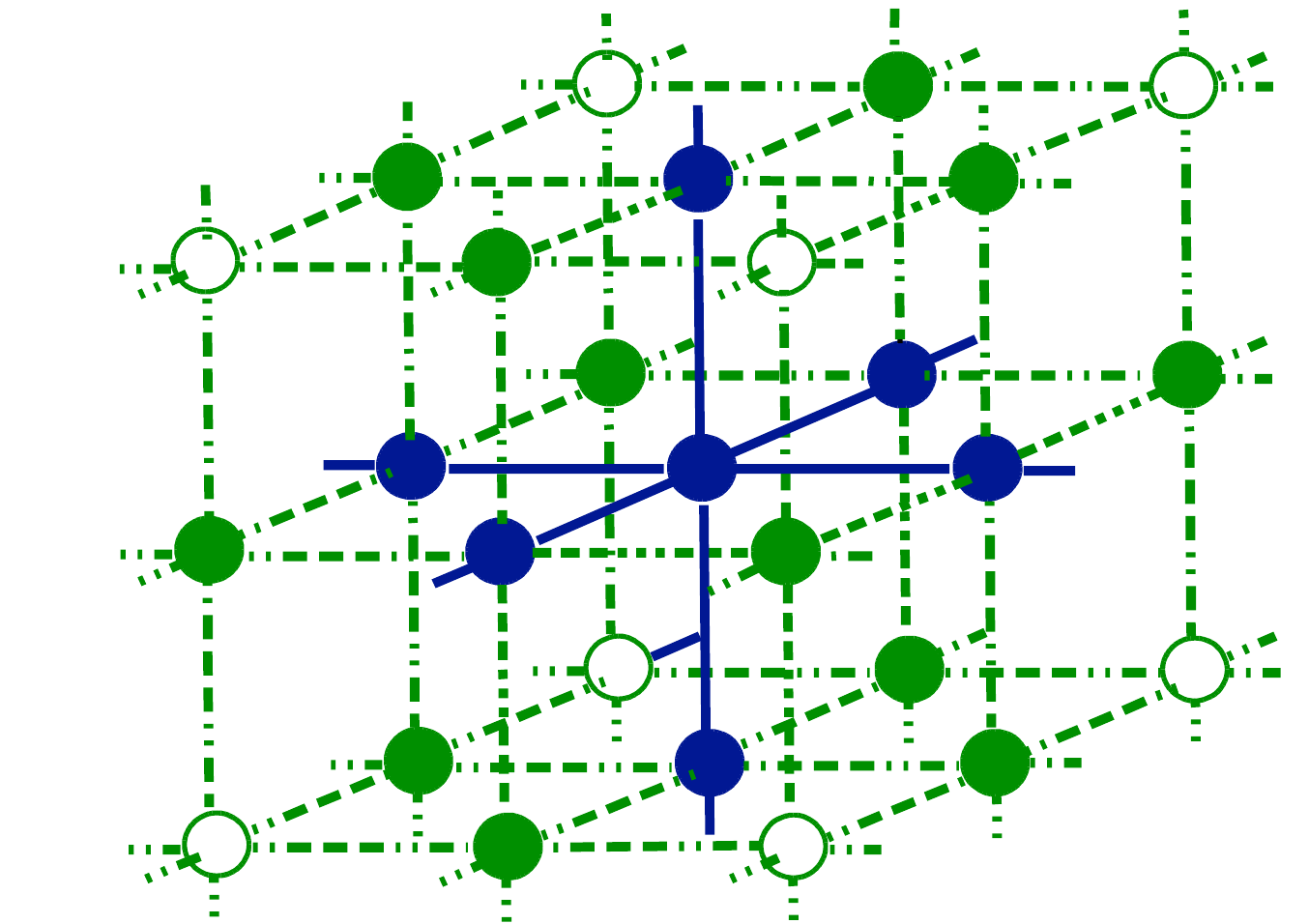}}
\caption{oscillations in the infinite weak component}\label{doublep1bis}
   \end{figure}
   
\begin{example}\rm
We finally consider a three-dimensional two-periodic geometry, with one strong connected component
pictured in Fig.~\ref{doublep1bis}.
Even in the absence of the forcing term $g$ we may have several competing microstructures
in the determination of $\varphi$. In Fig~\ref{doublep1bis} we have represented the uniform
data $u=+1$ on the strong component with solid circles, and a system of ferromagnetic connections
between strong and weak sites (positive coefficients) and of antiferromagnetic connections
between weak sites (a negative coefficient $\alpha$). Correspondingly, the minimal states 
have the value $+1$ on weak sites connected with the strong component (represented by
solid circles), and the value $-1$ on the other sites (represented by white circles).
Note that in this case the contribution of the weak phase is a constant.
\end{example}

\end{document}